\documentclass{article}

\input{amssym.def}
\input{amssym.tex}
\input{epsf.tex}
\input{xy}
\xyoption{all}

\usepackage{amsthm}
\usepackage{epsfig}
\usepackage{color}

\title{Commensurators of finitely generated non-free Kleinian groups}
\author{C. J. Leininger, D. D. Long \& A. W. Reid\thanks{All three
authors supported in part by the NSF}}

\def\H{\rm{H}}
\def\PSL{\rm{PSL}}
\def\Stab{\rm{Stab}}

\def\qed{ $\sqcup\!\!\!\!\sqcap$}

\def\sg{\mbox{\rm{sg}}}

\def\Area{\mbox{\rm{Area}}}

\newtheorem{theorem}{Theorem}[section]
\newtheorem{proposition}[theorem]{Proposition}
\newtheorem{lemma}[theorem]{Lemma}
\newtheorem{corollary}[theorem]{Corollary}

\begin{document}
\maketitle

%%%%%%%%%%%%%%%%%%%%%%%%%%%%%%%%%%%%%%%%%%%%%%%%%%%%%
%
%
% Introduction
%
%
%%%%%%%%%%%%%%%%%%%%%%%%%%%%%%%%%%%%%%%%%%%%%%%%%%%%%
\section{Introduction}

Let $G$ be a group and $\Gamma_1, \Gamma_2 <G$.  $\Gamma_1$ and $\Gamma_2$
are called commensurable if $\Gamma_1\cap \Gamma_2$
has finite index in both $\Gamma_1$ and $\Gamma_2$.
The {\em Commensurator} of a subgroup $\Gamma < G$
is defined to be:

$$C_G(\Gamma) = \{g\in G:g\Gamma g^{-1}~\hbox{is commensurable with}~\Gamma\}.$$

When $G$ is a semi-simple Lie group, and $\Gamma$ a lattice,
a fundamental dichotomy established by Margulis \cite{Ma}, determines
that $C_G(\Gamma)$ is dense in $G$ if and only if $\Gamma$ is arithmetic, and
moreover,
when $\Gamma$ is non-arithmetic, $C_G(\Gamma)$ is again a lattice.

Historically, the prominence of the commensurator was due in large part to
the density of the commensurator in the arithmetic setting being closely
related to the abundance of Hecke operators attached to arithmetic lattices.
These operators are fundamental objects
in the theory of automorphic forms associated to arithmetic lattices
(see \cite{Shi} for example).
More recently, the commensurator of various classes
of groups has come to the fore due
to its growing role in geometry, topology and geometric group
theory; for example in classifying lattices up to quasi-isometry,
classifying graph manifolds up to quasi-isometry, and understanding Riemannian
metrics admitting many ``hidden symmetries''
(for more on these and other topics see \cite{BB}, \cite{BeN}, \cite{FW1},
\cite{FW2}, \cite{LM}, \cite{Sch} and \cite{Sh}).

In this article, we will study $C_G(\Gamma)$ when $G=\PSL(2,{\bf C})$ and
$\Gamma$ a finitely generated non-elementary
Kleinian group.
In this setting, we will abbreviate the notation for the commensurator of
$\Gamma$ to $C(\Gamma)$.
In the case that $\Gamma$ is of finite co-volume and
non-arithmetic, identifying $C(\Gamma)$ has attracted considerable
attention (see for example \cite{GGH}, \cite{NR} and \cite{RW}).
Our focus here are those Kleinian groups $\Gamma$ for which
${\bf H}^3/\Gamma$
has infinite volume.
Henceforth, unless otherwise stated, the Kleinian groups $\Gamma$
that we consider will {\bf always have infinite co-volume}.

If $\Gamma$ is a Kleinian group, we denote by $\Lambda_\Gamma$ and
$\Omega_\Gamma$ the limit set and domain of discontinuity of $\Gamma$.
The Kleinian group $\Gamma$ is said to be of the {\em first kind}
(resp. {\em second kind}) if $\Omega_\Gamma=\emptyset$
(resp. $\Omega_\Gamma\neq \emptyset$).  The only known result for
Kleinian groups $\Gamma$ as above is due to L. Greenberg who proved
the following result (see \cite{Gre1} and \cite{Gre2}, and see \S 3.1
for a new proof when there are no parabolics).

\begin{theorem}
\label{greenberg}
Let $\Gamma$ be a finitely generated
non-elementary Kleinian group of the second kind and assume that $\Lambda_\Gamma$
is not a round circle. Then $[C(\Gamma):\Gamma]<\infty$. In particular,
$C(\Gamma)$ is a discrete
subgroup of $\PSL(2,{\bf C})$.\end{theorem}

The exclusion of the limit set being a round circle is to rule out
arithmetic Fuchsian subgroups of $\PSL(2,{\bf R})$, which by \cite{Ma}
have non-discrete commensurator.  Thus it remains to understand
$C(\Gamma)$ when $\Gamma$ is of the first kind.  These have been the
most difficult Kleinian groups to understand, even in the case that
$\Gamma$ is isomorphic to a closed surface group. Indeed, it is only
very recently that their geometry has been clarified; see \cite{Mi2},
\cite{BCM}, \cite{MM1}, \cite{A}, \cite{CG}.

Our main result is:

\begin{theorem}
\label{commensurator}
Let $\Gamma$ be a finitely generated torsion-free Kleinian
group of the first kind which is not a lattice. If $\Gamma$ is not free and contains no parabolic elements,
then $C(\Gamma)$ is discrete.

Furthermore $[C(\Gamma):\Gamma]=\infty$ if and only if
$\Gamma$ is a fibre group. In this case, $C(\Gamma)$ is a lattice.\end{theorem}

We recall that a  Kleinian group $\Gamma$ is  called a {\em
fibre group} if there is a finite volume hyperbolic 3-manifold
$M={\bf H}^3/\Delta$ which fibres over a $1$--orbifold such that $\Gamma < \Delta$ is the subgroup corresponding to a fibre.

The torsion-free assumption is simply to make the statement succinct. Since
every finitely generated Kleinian group contains a torsion free subgroup
of finite index, it is clear that Theorem \ref{commensurator} holds
in the presence of torsion with only a mild change of terminology for fibre groups in this setting.

We suspect that Theorem \ref{commensurator} holds for any finitely generated Kleinian group which is not a lattice, that is, without the hypothesis on parabolics or freeness.
In fact, the proof can be slightly modified to work without the assumption on parabolics, provided we further assume that $\Gamma$ is not the fundamental group of a compression body with toroidal lower boundary; see Section \ref{S:parabolics}.  However, there remain issues with generalizing the proof to the case of a free Kleinian group, with or without parabolics.

The proof of Theorem \ref{commensurator} makes
use of recent progress in understanding finitely generated geometrically
infinite Kleinian groups $\Gamma$. In particular, our proof appeals directly to
work of M. Mj on the
existence of ``Cannon-Thurston maps'' (see \cite{MM1} and \cite{MM2} and
\S 2 for more discussion on this), which in turn relies on existence of
``models for simply degenerate ends'' following Minsky, and
Brock, Canary and Minsky (see \cite{Mi2}, \cite{BCM} and \cite{BCM2}).
In turn, the complete classification of Kleinian groups via their
end invariants also uses the solution to the Tameness Conjecture
(\cite{A} and \cite{CG}), although in our case, the application
of the Cannon-Thurston map arises for geometrically infinite Kleinian
surface groups, and tameness of the quotients of ${\bf H}^3$ was already established
by work of Bonahon \cite{Bon}.

An important ingredient of our proof is of some independent
interest. We defer a precise statement to Theorem \ref{notconical},
but it can be described informally as follows.  The Cannon-Thurston map
$\pi : S^1 \longrightarrow S^2$ is obtained via a certain
decomposition map; in particular it picks out in the image a
collection of exceptional points, namely those points $p \in S^2$ for
which $|\pi^{-1}(p)| > 1$.  The result of \ref{notconical} is:

\begin{theorem}
\label{notconicalIntro}
The exceptional points of the Cannon-Thurston map are not conical limit points.
\end{theorem}
By the solution to the Tameness Conjecture (\cite{A} and \cite{CG}),
and the Covering Theorem of Canary \cite{Ca1} (see Theorem \ref{canary}),
every finitely generated geometrically infinite subgroup of a uniform
lattice in $\PSL(2,{\bf C})$ is isomorphic to a closed surface group. Thus
a corollary of Theorems \ref{greenberg}
and \ref{commensurator} is the following result. This extends a result
of the third author (see \cite{Ca2} Theorem 8.7).

\begin{corollary}
\label{geominf}
Let $\Delta$ be a uniform lattice in $\PSL(2,{\bf C})$, and
$\Gamma < \Delta $ a finitely generated subgroup of infinite index
for which $\Lambda_\Gamma$ is not a circle.
Then $C(\Gamma)$ is a discrete subgroup of $\PSL(2,{\bf C})$.\end{corollary}

\noindent{\bf Acknowledgements:}~{\em The authors thank Yehuda Shalom for
useful correspondence on the subject of this paper, and Dick Canary, Yair Minsky and Hossein Namazi
for helpful conversations.}

%%%%%%%%%%%%%%%%%%%%%%%%%%%%%%%%%%%%%%%%%%%%%%%%%%%%%%%%%%%%%%%%%%%%%%%%%%%%%%%%%%%%%
%%%%%%%%%%%%%%%%%%%%%%%%%%%%%%%%%%%%%%%%%%%%%%%%%%%%%%%%%%%%%%%%%%%%%%%%%%%%%%%%%%%%%
\section{Ending Laminations and Cannon-Thurston Maps} \label{S:background}
%%%%%%%%%%%%%%%%%%%%%%%%%%%%%%%%%%%%%%%%%%%%%%%%%%%%%%%%%%%%%%%%%%%%%%%%%%%%%%%%%%%%%
%%%%%%%%%%%%%%%%%%%%%%%%%%%%%%%%%%%%%%%%%%%%%%%%%%%%%%%%%%%%%%%%%%%%%%%%%%%%%%%%%%%%%

The proof of Theorem \ref{commensurator} requires some background
about hyperbolic 3-manifolds
and from the theory of laminations, ending laminations and Cannon-Thurston
maps. We summarize what we need here.

%%%%%%%%%%%%%%%%%%%%%%%%%%%%%%%%%%%
\subsection{} \label{S:laminations}
%%%%%%%%%%%%%%%%%%%%%%%%%%%%%%%%%%%

We begin by recalling some basic structure for laminations
on a closed surface $S$ (see \cite{CB} for further details).

Let $S$ be a closed orientable surface of genus $g\geq 2$ equipped
with a fixed complete hyperbolic metric of constant curvature $-1$. A
(geodesic) lamination $L$ on $S$ is a closed subset foliated by geodesics.
A component of $S\setminus L$ is called a {\em principal region} for $L$, of
which there are only finitely many (see \cite[Lemma 4.3]{CB}).
A lamination $L \subset S$
is {\em filling} if it has no proper sublaminations and all principal
regions are ideal polygons.
In particular $L$ has no closed leaves.
%A {\em boundary leaf} of a principal region $P$ for $L$ is a leaf $\ell$ of
%$L$ such that for all $x\in \ell$ there exists an $\epsilon>0$ such that
%$B_\epsilon(x)\cap P$ contains a component $B$ of
%$B_\epsilon(x)\setminus \sigma$, where $\sigma$ is the length $2\epsilon$
%segment of $\ell$ centered at $x$, for which the closure $\overline{B}$
%contains $\sigma$.

Given a lamination $L \subset S$, we let $\widetilde L \subset
\widetilde S$ denote the preimage of $L$ in the universal cover
$\widetilde S$ of $S$.  If $\widetilde \ell$ (respectively $\widetilde
P$) is a leaf (respectively, principal region) of $\widetilde L$, then
we write $\partial_\infty \widetilde \ell$ (respectively,
$\partial_\infty \widetilde P$) for the intersection of the closure of
$\widetilde \ell$ (respectively, $\widetilde P$) with $S^1_\infty$,
the circle at infinity of $\widetilde S$.

%%%%%%%%%%%%%%%%%%%%%%%%%%%%%%
\subsection{} \label{S:ending}
%%%%%%%%%%%%%%%%%%%%%%%%%%%%%%

Let $\Gamma$ be a finitely generated non-elementary Kleinian group
without parabolics and $N={\bf H}^3/\Gamma$. Throughout the paper the
{\em convex core} of $N$ will be denoted by $C_N$; i.e. $C_N =
CH(\Lambda_\Gamma)/\Gamma$, where $CH(\Lambda_\Gamma)$ is the convex
hull of $\Lambda_\Gamma$ in ${\bf H}^3$. If $C_N$ is compact, then
$\Gamma$ is said to be {\em geometrically finite} (or convex
cocompact), otherwise $\Gamma$ is called geometrically infinite.  We
will also write $M \subset N$ for a fixed {\em compact core} of $N$.
This is a compact submanifold of $N$ for which the inclusion map is a
homotopy equivalence (see Scott \cite{Sc1}).

Given a component $S \subset \partial M$ let $U \subset (N\setminus M)$ be
the component with $S \subset \overline U$.
This $U$ is a neighborhood of an end $E$ of $N$, and we say that $E$ {\em abuts} $S$.
This defines a bijection between the set of ends of $N$ and the components of $\partial M$.
According to the Tameness Theorem (\cite{A} and \cite{CG}) $U \cong S \times (0,\infty)$.

Following Thurston \cite{Th}, the end $E$ of $N$ abutting $S \subset \partial M$
is called {\em geometrically finite} if $E$ contains a neighbourhood $U$ such that $U \cap C_N = \emptyset$.
The group $\Gamma$ is geometrically finite if and only if all ends of $C_N$ are geometrically finite.  In this case, we can take $M = C_N$.

An end $E$ of $N$ is {\em simply degenerate} if there exists a sequence of simple closed curves $\{\alpha_i\}$ in $S$
such that the geodesic representatives $\{\alpha_i^*\}$ in $N$ exit the end $E$.
The Tameness Theorem implies that every end is either geometrically finite or simply degenerate.

Furthermore, Thurston \cite{Th2} and Canary \cite{Ca0} showed how
to associate to a simply degenerate end $E$, a
lamination (the so-called {\em ending lamination} of E)
which is a limit of the simple closed curves $\{\alpha_i\}$
(in an appropriate topology).  We shall denote the ending lamination
associated to the simply degenerate end $E$ by $\nu_E$.
If $\nu_E$ is an ending lamination then it
is known to be filling (see \S 8 of \cite{Ca0} for example).

%%%%%%%%%%%%%%%%%%%%%%%%%%%%%%%%
\subsection{} \label{S:covering}
%%%%%%%%%%%%%%%%%%%%%%%%%%%%%%%%

An important tool for us is Canary's Covering Theorem \cite{Ca1},\cite{A} (which extended a result of
Thurston \cite{Th}).  The version we state here can be found in \cite{CL}.
%%%%%%%%%%%%%%%%%%%%%%%%%%%%%%%%%%%%%%%
%
%
% Canary covering theorem
%
%
%%%%%%%%%%%%%%%%%%%%%%%%%%%%%%%%%%%%%%%
\begin{theorem}
\label{canary}
Suppose that $\hat N = {\bf H}^3/\Gamma$ is a hyperbolic 3-manifold with no cusps and finitely generated fundamental group and $N$
is a hyperbolic 3-orbifold. If $p:\widehat N \to N$ is a cover which is infinite-to-one on a neighborhood $U$ of
a simply degenerate end of $\hat N$, then $N$ is closed and has a finite manifold cover $N'={\bf H}^3/\Gamma' \to N$
such that either
\begin{enumerate}
\item
$N'$ fibers over the circle and $\hat N$ is the cover associated
to a fiber subgroup of $\Gamma'$, or
\item
$N'$ fibers over the orbifold
$S^{1}/\langle z \mapsto \overline{z} \rangle$ and
$\hat N$ is the cover of $N'$ associated to a singular fiber subgroup
of $\pi_1(N')$.
\end{enumerate}
\end{theorem}
In the conclusion of the theorem, the group $\Gamma$ is a fibre group and $N$ is an ${\bf R}$--bundle over a surface.  In the first case, $N$ is a product, while in the second, it is a twisted ${\bf R}$--bundle.

%%%%%%%%%%%%%%%%%%%%%%%%%%
\subsection{} \label{S:CT}
%%%%%%%%%%%%%%%%%%%%%%%%%%

Now suppose that $\Gamma$ is a {\em surface group}, so that the compact core is a product $M = S \times [0,1]$, where $S$ is closed orientable surface of genus at least $2$.   In this case, there are exactly two ends $E_+$ and $E_-$.  If $\Gamma$ is geometrically infinite, one or both of the ends is simply degenerate and we say that $\Gamma$ is {\em singly degenerate} or {\em doubly degenerate} in these two cases, respectively.
For the remainder of this section, we assume that $\Gamma$ is doubly degenerate.  For example, $\Gamma$ may be a fibre group of the first type described in the previous subsection.
The ending laminations $\nu_+ = \nu_{E_+}$ and $\nu_- = \nu_{E_-}$ can be viewed as laminations on $S$.

If $\Gamma$ contains a finite index subgroup which is the fundamental group of an orientable surface, then we will say that $\Gamma$ is virtually a surface group. It follows from \cite[Theorem 10.5]{He} that if $\Gamma$ is virtually a surface group then $N = {\bf H}^3/\Gamma$ is an ${\bf R}$--bundle over a surface, and hence either $\Gamma$ or a canonical index two subgroup $\Gamma_0 < \Gamma$ is a surface group.

The inclusion map $S \hookrightarrow S \times \{1/2\} \subset M \subset N$
induces a representation
\[\rho:\pi_1(S) \to \Gamma < \PSL(2,{\bf C})\]
and lifts to an equivariant map
\[\iota:\widetilde S \to {\bf H}^3.\]
The equivariance is with respect to $\rho$ via the action of $\pi_1(S)$ on $\widetilde S$ by covering transformations and $\Gamma$ on ${\bf H}^3$:
$$\iota \circ \gamma = \rho(\gamma) \circ \iota,~ \hbox{for all}~\gamma\in \pi_1(S).$$

In \cite{CT}, Cannon and Thurston proved that if $\Gamma$ is a fibre group, then $\iota$ admits a continuous equivariant extension to the compactifications
\[\hat \iota: \widetilde S \cup S^1_\infty \rightarrow {\bf H}^3 \cup S^2_\infty.\]
The existence of such an extension was subsequently proven by Minsky
\cite{Mi} replacing the fibre assumption with the weaker assumption of
{\em bounded geometry} (i.e. a there is a lower bound to the injectivity
radius).  Existence in the general case was proven more
recently by Mj \cite{MM1}, and $\hat \iota$ is referred to as a {\em
  Cannon--Thurston map} for $\Gamma$ or $\rho$.  We state Mj's Theorem
here together with the description of the map that we will need (see
\cite[Theorem 1.3]{MM2}).

\begin{theorem}
\label{CTMapsforsurfaces}
Given a representation $\rho:\pi_1(S) \to  \PSL(2,{\bf C})$ with doubly degenerate image $\Gamma = \rho(\pi_1(S))$ as above, there exists a Cannon--Thurston map
\[\hat \iota:\widetilde S \cup S^1_\infty \to {\bf H}^3 \cup S^2_\infty.\]
If $\nu_\pm$ are the ending laminations and $a,b \in S^1_\infty$, then $\hat \iota(a) = \hat \iota(b)$ if and only if $a$ and $b$ are either
ideal end points of a leaf, or ideal endpoints of a principal region of one of $\widetilde \nu_+$ or $\widetilde \nu_-$.
\end{theorem}

It is straightforward to show that if $a$ and $b$ are either
ideal end points of a leaf of $\widetilde \nu_\pm$ or ideal endpoints of a principal
region of $\widetilde \nu_\pm$ then they are identified by the Cannon-Thurston map (see
Lemma 3.5 of \cite{MM0}).  The hard part is to show that this is all that is collapsed.

\section{Proof of Theorem \ref{commensurator}}

The proof of Theorem \ref{commensurator} proceeds by examining two
cases. The first is the case where neither $\Gamma$ nor an index two subgroup is a surface group.

%%%%%%%%%%%%%%%%%%%%%%%%%%%%%%%%%%%%%%%%%%%%%%%%%%%%%%%%%%%%%%%%%%%%%%%%%%%%%%%%%%%%%%%%%%
\subsection{The case that $\Gamma$ is not virtually a surface group.} \label{S:not surface}
%%%%%%%%%%%%%%%%%%%%%%%%%%%%%%%%%%%%%%%%%%%%%%%%%%%%%%%%%%%%%%%%%%%%%%%%%%%%%%%%%%%%%%%%%%

We will prove the following result.  This implies the aforementioned
case of Theorem \ref{commensurator}, and in addition gives a new proof
Greenberg's result when there are no parabolics (cf. Theorem \ref{greenberg}).

%%%%%%%%%%%%%%%%%%%%%%%%%%%%%%%%%%%
%
% The Temptation theorem
%
%%%%%%%%%%%%%%%%%%%%%%%%%%%%%%%%%%%
\begin{theorem}
\label{notsurface}
Let $\Gamma$ be a finitely generated torsion-free Kleinian group without
parabolic elements for which $\Lambda_\Gamma$ is not contained in a round circle.  In addition, assume
that if $\Gamma$ is of the first kind it is not a lattice or isomorphic to a free group or a surface group.

Then $C(\Gamma)$ is discrete and $[C(\Gamma):\Gamma] < \infty$.
\end{theorem}
The hypothesis that $\Lambda_\Gamma$ is not contained in a round circle can be replaced by the assumption that $\Lambda_\Gamma$ is not equal to a round circle.  The proof is a (simpler) two-dimensional version of the argument we give here, and we leave the details to the interested reader.\\

\noindent{\bf Proof:}~Let $q : {\bf H}^3 \longrightarrow N$ denote the
universal covering map.

Suppose first that $\Gamma$ is of the second kind so that
$\Omega_\Gamma \neq \emptyset$.  In this case, there is at least one
geometrically finite end of $N={\bf H}^3/\Gamma$.  Thurston \cite{Th}
showed that each component of $\partial C_N$ is a pleated surface (see
also \cite{EM}) which we may assume is also a component of $\partial
M$.  For each geometrically finite end $E$ of $N$, we let ${\cal P}_E$
denote this pleated surface which $E$ abuts.  Let $E_1,\ldots,E_k$
denote the geometrically finite ends of $N$ and write
\[ {\cal P}_N = \{ {\cal P}_{E_1} , \ldots , {\cal P}_{E_k}\}. \]

The preimage of this finite set of pleated surfaces in ${\bf H}^3$ is a
locally finite collection of connected pleated surfaces in ${\bf H}^3$
\[ {\cal P}_\Gamma = q^{-1}({\cal P}_N).\]
If $\hat \Gamma < \Gamma$ is a finite index subgroup then
\begin{equation} \label{fi eqn1}
{\cal P}_\Gamma = {\cal P}_{\hat \Gamma}
\end{equation}
To see this, suppose that $\hat E_1,\ldots,\hat E_r$ are the ends of
$\hat N$ with neighborhoods covering a neighborhood of $E$ via the
associated covering $p:\hat N \to N$.  Then
\[p^{-1}({\cal P}_E) = \{{\cal P}_{\hat E_1}, \ldots, {\cal P}_{\hat E_r} \}\]
From this it follows that $p^{-1}({\cal P}_N) = {\cal P}_{\hat N}$, and consequently (\ref{fi eqn1}) holds.
(Alternatively, we can simply observe that ${\cal P}_\Gamma = {\cal P}_{\hat \Gamma}$ is the set of boundary components of $CH(\Lambda_\Gamma) = CH(\Lambda_{\hat \Gamma})$, however the first argument illustrates the idea needed for the case $\Omega_\Gamma = \emptyset$ below).

Given $\gamma \in \PSL(2,{\bf C})$ observe that
\begin{equation}
\label{conj eqn1}
{\cal P}_{\gamma \Gamma \gamma^{-1}} = \gamma({\cal P}_\Gamma).
\end{equation}
Furthermore, if $\gamma \in C(\Gamma)$, $\hat \Gamma = \gamma \Gamma \gamma^{-1} \cap \Gamma$ has finite index in both $\Gamma$ and $\gamma \Gamma \gamma^{-1}$, so combining (\ref{fi eqn1}) and (\ref{conj eqn1}), we have
\[ \gamma({\cal P}_\Gamma) = {\cal P}_\Gamma.\]
Therefore,
\[ C(\Gamma) < \Stab({\cal P}_\Gamma) = \{g \in \PSL(2,{\bf C}) \, | \, g {\cal P}_\Gamma = {\cal P}_\Gamma\}, \]
and so discreteness of $C(\Gamma)$ follows from the next claim.\\

\noindent {\bf Claim:}~$\Stab({\cal P}_\Gamma)$ is discrete.\\
% \marginpar{\tiny I added these two paragraphs to indicate where we
%   are using the fact that $\Lambda_\Gamma$ is not contained in a
%   round circle... we can remove it if it's obvious.}

\noindent {\bf Proof of Claim:}
Each surface in ${\cal P}_\Gamma$, contains countably many flat pieces, each of which is contained in a hyperbolic plane ${\bf H}^2 \subset {\bf H}^3$.  Let $\{ {\bf H}^2_i\}$ be the set of all such hyperbolic planes with ${\bf H}_i^2$ containing the flat piece $F_i$ of some pleated surface in ${\cal P}_\Gamma$.  The intersection $\cap_i {\bf H}^2_i$ is invariant by $\Stab({\cal P}_\Gamma)$ (and so also $\Gamma$), and since $\Lambda_\Gamma$ is not contained in a round circle, it follows that this intersection is empty.   Therefore, there are a finite set of flat pieces, say $F_1,\ldots,F_k$, of surfaces in ${\cal P}_\Gamma$ so that the associated hyperbolic planes intersect trivially
\[ {\bf H}_1^2 \cap \cdots \cap {\bf H}_k^2 = \emptyset.\]

Now suppose that $\{\gamma_n\}_{n=1}^\infty \subset \Stab({\cal P}_\Gamma)$ is a sequence converging to the identity.  Since ${\cal P}_\Gamma$ is locally finite, after passing to a subsequence we can assume that each $F_i$ is invariant by $\gamma_n$ for $i = 1,\ldots,k$ and all $n$.  Therefore, ${\bf H}^2_i$ is invariant by $\gamma_n$ for each $i =1,\ldots,k$ and all $n$.  Since ${\bf H}^2_1 \cap \cdots \cap {\bf H}^2_k = \emptyset$ it follows that $\gamma_n$ is the identity, and hence $\Stab({\cal P}_\Gamma)$ is discrete.\qed \\

We now assume that $\Gamma$ is of the first kind.  The idea of the proof is similar to the previous case: to each end $E$ of $N$ we associate a finite, nonempty set of pleated surfaces ${\cal P}_E$ in $N$.  This will be done so that the preimage ${\cal P}_\Gamma = q^{-1}({\cal P}_N)$ satisfies (\ref{fi eqn1}) and (\ref{conj eqn1}).  The proof can then be completed exactly as above.

Given an end $E \subset N$, let $S \subset \partial M$ be the component which it abuts.
According to \cite[Theorem 2.1]{Bon0} there is a unique (up to isotopy) compression body $B_S \subset M$ such that homomorphism $\pi_1(B_S) \to \pi_1(M) = \Gamma$ induced by inclusion is injective.  The lower boundary $\partial_- B_S$ is a finite disjoint union of surfaces, which is nonempty since $\Gamma$ is not free.  We let $\Gamma_S^1,\ldots,\Gamma_S^k < \Gamma$ denote the fundamental groups of the components---these groups inject into $\pi_1(B_S)$ and so also into $\Gamma$.  Let $N_S^i = {\bf H}^3/\Gamma_S^i \to N$ be the associated covers.

If $S$ is incompressible, then $B_S \cong S \times [0,1]$, $\partial_- B_S$ has one component and the cover $N_S^1 \to N$ is the cover corresponding to $\pi_1(S) < \Gamma$.  Since $\Gamma$ is not virtually a surface group, this is an infinite sheeted cover.  By Theorem \ref{canary}, $N_S^1$ has at least one geometrically finite end (in fact, it has exactly one geometrically finite end, and one simply degenerate end with a neighborhood that maps isometrically onto a neighborhood of $E$ in $N$).
If $S$ is compressible, then each $\Gamma_S^i$ has infinite index in $\pi_1(B_S)$, and hence also in $\Gamma$.  Appealing to Theorem \ref{canary} again we see that $N_S^i$ has at least one geometrically finite end for each $1 \leq i \leq k$.

For each $i$, we obtain one or two pleated surfaces in $N$ as the image of $\partial C_{N_S^i}$ by the covering map $N_S^i \to N$.  Let ${\cal P}_E$ denote the set of all such, as $i$ ranges from $1$ to $k$.  If $E_1,\ldots,E_n$ are the ends of $N$, then set
\[ {\cal P}_N = \{{\cal P}_{E_1} , \ldots, {\cal P}_{E_n}\}\]
and put
\[ {\cal P}_\Gamma = q^{-1}({\cal P}_N).\]
With this construction of pleated surfaces, it is straightforward to check that (\ref{fi eqn1}) and (\ref{conj eqn1}) hold.  The remainder of the proof now proceeds as in case of $\Omega_\Gamma \neq \emptyset$, proving that $C(\Gamma)$ is discrete.\\

To prove the final statement that $[C(\Gamma):\Gamma] < \infty$, let ${\cal P}$ be the union of all the pleated surfaces in ${\cal P}_\Gamma$.  Then
\[ [C(\Gamma):\Gamma] = \frac{\Area({\cal P}/\Gamma)}{\Area({\cal P}/C(\Gamma))} < \infty.\]
This completes the proof. \qed.

%%%%%%%%%%%%%%%%%%%%%%%%%%%%%%%%%%%%%%%%%%%%%%%%%%%%%%%%%%%%%
\subsection{The surface group case.} \label{S:surface group}
%%%%%%%%%%%%%%%%%%%%%%%%%%%%%%%%%%%%%%%%%%%%%%%%%%%%%%%%%%%%%

Given Theorem \ref{greenberg} and Theorem \ref{notsurface}
the proof of Theorem \ref{commensurator} will be completed upon establishing:

%%%% we still need to show that

\begin{theorem}
\label{surfacegroup}
Let $S$ be a closed orientable surface,
$\rho:\pi_1(S)\rightarrow \PSL(2,{\bf C})$
a faithful discrete representation and
let $\Gamma  = \rho(\pi_1(S))$. Assume that $\Gamma$ is doubly degenerate.

Then $C(\Gamma)$ is discrete and finitely generated. Furthermore, $[C(\Gamma):\Gamma]=\infty$ if and only if
$\Gamma$ is a fibre group. In this case, $C(\Gamma)$ is a lattice.
\end{theorem}

Before proving Theorem \ref{surfacegroup} we recall some terminology
that will be used.\\[\baselineskip]
\noindent{\bf Definition:}~{\em Let $\Gamma$ be a Kleinian group. A point
$x\in \Lambda_\Gamma$ is called a conical limit point of $\Gamma$ if for some (and hence
every) geodesic ray $\tau$ in ${\bf H}^3$ ending at $x$, there is a compact set
$K\subset {\bf H}^3$ such that $\{\gamma\in \Gamma : \gamma(\tau)\cap K \neq \emptyset\}$
is infinite.}\\[\baselineskip]
We will make use of the following equivalent version of a conical limit
point. We sketch a proof for convenience.

\begin{lemma}
\label{conical}
Let $\Gamma$ be as above. The limit point $x$ is a conical limit
point for $\Gamma$ if and only
if $\Gamma$ contains a sequence of elements $\{\gamma_m\}$ such that
the following conditions hold:

\item (i) $\gamma_m(x)\rightarrow x_\infty^+ \in S_\infty^2$, and

\item (ii) $\gamma_m|_{S_\infty^2\setminus \{x\}}$ converges uniformly on compact
sets to a constant map with value $x_\infty^- \neq x_\infty^+$.\end{lemma}

\noindent{\bf Proof:}(Sketch)~It is convenient to
pass to the Ball Model ${\bf B}$ of ${\bf H}^3$. Then, as in
\cite[pp.~3--4]{BM},
that $x$ is a conical limit point
is equivalent to:\\[\baselineskip]
{\em For some other point $y\in S_\infty^2$,
there exists a sequence of distinct elements $\gamma_n \in \Gamma$
such that
$$|\gamma_n(y)-\gamma_n(x)|\geq \delta>0.$$}
Briefly, taking $\sigma=[y,x]$ to be a diameter of $\bf B$, the condition that
$|\gamma_n(y)-\gamma_n(x)|\geq \delta > 0$ can be
seen to be equivalent to the statement that the geodesics $\gamma_n(\sigma)$
remain close to $0$, and so meet some closed ball centered at $0$ infinitely
often.

Now as discussed in \cite[p.~122]{Mas}, this version of conical limit point
can be shown to be equivalent to the statement of Lemma \ref{conical}.
\qed\\[\baselineskip]
We now commence with the proof of Theorem \ref{surfacegroup}.\\[\baselineskip]
\noindent{\bf Proof of Theorem \ref{surfacegroup}:}~Let $N={\bf H}^3/\Gamma$
with ending laminations $\nu_+$ and $\nu_-$. These laminations have the
property that they {\em bind}
in the following sense (see \cite{Mi} and \cite{Ke}, although in
these papers, this property is called filling which we wish to
avoid having used it earlier):\\[\baselineskip]
\noindent {\em The complement $S\setminus (\nu_+\cup \nu_-)$ is a union
of disks, each with boundary a finite collection of compact arcs coming
alternately from $\nu_+$ and $\nu_-$.  All but a finite number of these
disks have four sides.}\\[\baselineskip]
We denote by $D_1,\ldots ,  D_m$ those finite number of disks
which do not have four sides.
Fix one of these disks $D$, which is contained in the intersection of two principal regions $P_+$ of $\nu_+$ and $P_-$ of $\nu_-$.
Let $\widetilde D$ denote some lift of $D$ to the universal cover $\widetilde S$ of $S$.  Associated to $\widetilde D$ is a unique pair of lifts $\widetilde P_+$ and $\widetilde P_-$ of $P_+$ and $P_-$, respectively, for which $\widetilde D = \widetilde P_+  \cap \widetilde P_-$.

Let $\pi$ denote the restriction of the
Cannon-Thurston map $\hat{\iota}$ to $S^1_\infty$.
Theorem \ref{CTMapsforsurfaces} shows that the map $\pi$ collapses
$\partial_\infty \widetilde P_+$ and $\partial_\infty \widetilde P_-$ to points $x_+$ and $x_-$ in
$\Lambda_\Gamma = S^2_\infty$. Since the laminations bind, the closure of $\widetilde{D}$ is
compact and therefore $\partial_\infty \widetilde P_+ \cap \partial_\infty \widetilde P_- = \emptyset$.
It follows that $x_+$ and $x_-$ are distinct points in $S^2_\infty$.
We call these points {\em special points}.

Doing this for each of $D_1,\ldots,  D_m$ determines a finite number
of pairs of special points
$$(x_+^{(1)},x_-^{(1)}), (x_+^{(2)},x_-^{(2)}),\ldots (x_+^{(m)},x_-^{(m)})$$
in $S^2_\infty$.

Connect $x_+^{(j)}$ and $x_-^{(j)}$ by a geodesic $\xi^{(j)}$ in ${\bf
  H}^3$ for each $j=1,\ldots ,m$ and let ${\cal X}(\Gamma)$ be the
$\Gamma$-orbit of this collection of geodesics.  We shall
refer to a geodesic in ${\cal X}(\Gamma)$ as a {\em special geodesic}.
We record the following observation that is crucial in what follows.
%%%%%%%%%%%%%%%%%%%%%%%%%%%%%%%
%
%
% Not conical limit points
%
%
%%%%%%%%%%%%%%%%%%%%%%%%%%%%%%%
\begin{theorem}
\label{notconical}
Let $x\in S^2_\infty$ have the property that $|\pi^{-1}(x)| > 1$. Then
$x$ is not a conical limit point. In particular,
for each $i=1,2, \ldots m$, $x_+^{(i)}$ and $x_-^{(i)}$ are not conical limit
points.\end{theorem}
\noindent{\bf Proof:}~Given $x$ as in the hypothesis of
Theorem \ref{notconical}, Theorem \ref{CTMapsforsurfaces}
shows that there are $a_1,\ldots,a_k \in S^1_\infty$ ($k\geq 2$),
which are endpoints of a leaf or
ideal vertices of a principal region of $\widetilde \nu_+$ or $\widetilde \nu_-$ and which satisfy

$$x = \pi(a_1) = \ldots = \pi(a_k).$$

Moreover, by Theorem \ref{CTMapsforsurfaces}, these are the
only identifications.
We also recall that by definition of a Cannon--Thurston map, $\pi$ is equivariant with respect to $\rho$ via the action of $\pi_1(S)$ on $S^1 = S^1_\infty$ and $\Gamma$ on $S^2 = S^2_\infty$.

We now argue by contradiction, and assume that $x$ is a conical limit point.
Lemma \ref{conical} provides a sequence of distinct elements $\{\gamma_n\}_{n=1}^\infty \subset \Gamma$ so that
\begin{enumerate}
\item {\em $\gamma_n(x) \to x_\infty^+ \in S^2$ as $n \to \infty$, and}
\item {\em $\gamma_n|_{S^2 - \{x\}}$ converges uniformly on compact sets to a constant map with value $x_\infty^- \neq x_\infty^+$.}
\end{enumerate}

Let $\{g_n\}_{n=1}^\infty \subset \pi_1(S)$ be such that
$\rho(g_n) = \gamma_n$, and pass to a subsequence if necessary so that
as $n \to \infty$, $g_n (a_1) \to a_\infty \in S^1$ for some
$a_\infty \in S^1$.  Since $\pi$ is continuous and equivariant with
respect to $\rho$, we have
\[ \pi(a_\infty) = \lim_{n \to \infty} \pi(g_n(a_1)) = \lim_{n \to \infty} \rho(g_n)(\pi(a_1)) = \lim_{n \to \infty} \gamma_n(x) = x_\infty^+.\]

Next, let $b,c \in S^1 \setminus \{a_1,\ldots,a_k\}$ be distinct points
bounding an interval $I \subset S^1$ containing $a_2,\ldots,a_k$ and
\textit{not} containing $a_1$.  Pass to a further subsequence if
necessary so that as $n \to \infty$, $g_n(b) \to b_\infty$ and
$g_n(c) \to c_\infty$ for some points $b_\infty,c_\infty \in
S^1$.
\begin{figure}[htb]
\centerline{}
\centerline{}
\begin{center}
%\ \psfig{file=circle.eps,height=4cm} \caption{}
\includegraphics[height=4cm]{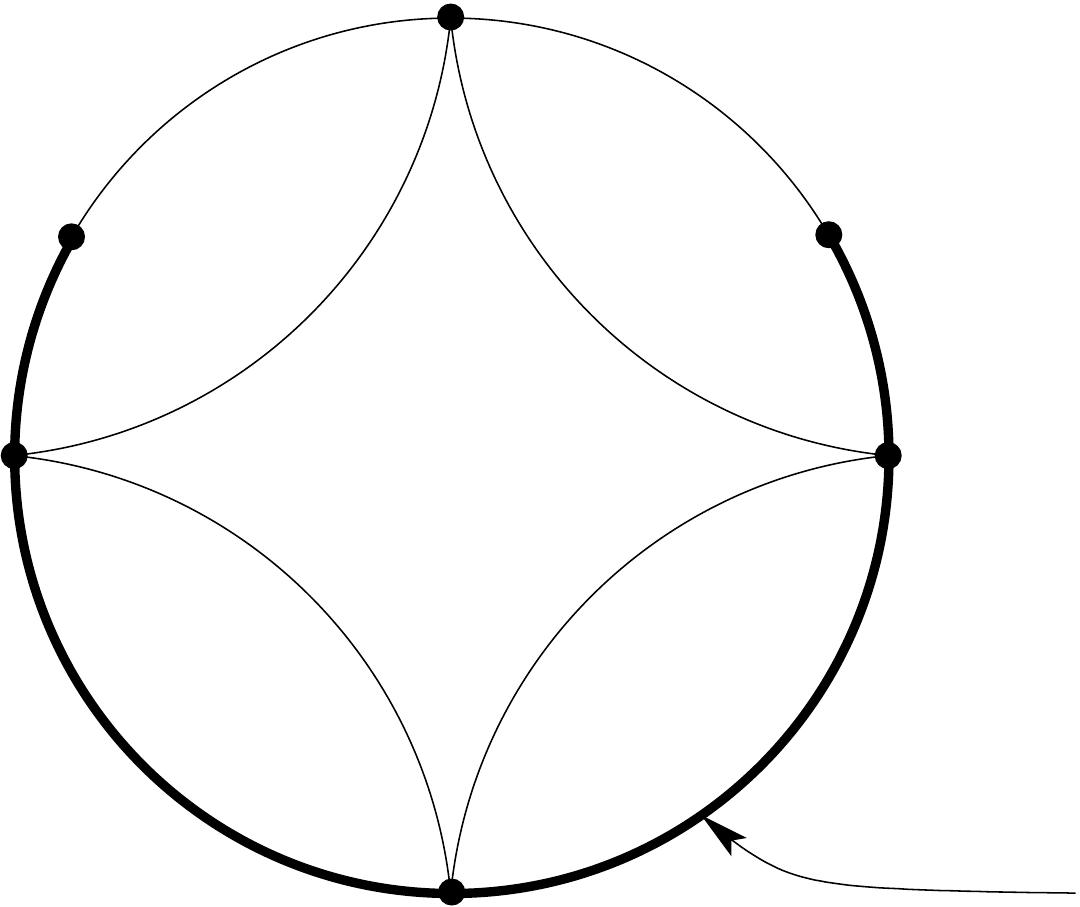}\caption{}
\label{circle}
\end{center}
  \setlength{\unitlength}{1cm}
  \begin{picture}(0,0)(0,0)
    \put(7.1,5.8){$a_1$}
    \put(4.8,3.5){$a_2$}
    \put(7.1,1.25){$a_3$}
    \put(9.4,3.5){$a_4$}
    \put(5.35,4.7){$b$}
    \put(9.03,4.67){$c$}
    \put(10.2,1.5){$I$}
  \end{picture}
\end{figure}

Appealing to continuity and equivariance of $\pi$ again we have
\begin{equation} \label{eq1}
\pi(b_\infty) = \lim_{n \to \infty} \pi (g_n (b)) = \lim_{n \to \infty} \gamma_n (\pi (b)) = x_\infty^- \neq x_\infty^+ = \pi(a_\infty)
\end{equation}
and
\begin{equation} \label{eq2}
\pi(c_\infty) = \lim_{n \to \infty} \pi (g_n (c)) = \lim_{n \to \infty} \gamma_n (\pi (c)) = x^-_\infty \neq x^+_\infty = \pi(a_\infty).
\end{equation}

Since the action of $\pi_1(S)$ on $\widetilde S \cup S^1$ is Fuchsian, a well-known characterization of Fuchsian
groups (see \cite[p.~484]{Sc})
%
%
%
% that follows from (\cite{Ga} and \cite{CJ})
%
%
%
now allows us to deduce that $\pi_1(S)$ acts properly discontinuously on
the space $(S^1\times S^1 \times S^1)\setminus \Delta$ where

$$\Delta=\{(p,q,r)\in S^1\times S^1 \times S^1 :~~\hbox{two of }p, q,
                                     r\hbox{ are equal}\}.$$

\noindent Hence $(a_\infty, b_\infty, c_\infty)$
cannot be a triple of distinct points in $S^1$ since it is the limit
of $\{g_n(a_1,b,c)\}_{n=1}^\infty$.  From (\ref{eq1}) and
(\ref{eq2}), $a_\infty \neq b_\infty$ and $a_\infty \neq c_\infty$.
Therefore, $b_\infty = c_\infty$.

It follows that after passing to another subsequence if necessary,
$\{g_n\}$ restricted to one of the intervals bounded by $b$ and
$c$, call it $J$, must converge to a constant map with value $b_\infty
= c_\infty$ (this is because each $g_n$ is a homeomorphism, and
the images of the endpoints are converging to the same point).  We
cannot have $a_1 \in J$ since $g_n(a_1) \to a_\infty \neq
b_\infty$ as $n \to \infty$.  Therefore, $J = I$, which is to say,
$g_n$ restricted to $I$ converges to the constant $b_\infty$.
However, this interval contains $a_2$, and so $g_n(a_2) \to
b_\infty$ as $n \to \infty$.  Another application of continuity and
equivariance of $\pi$ tells us
\[ \pi(b_\infty) = \lim_{n \to \infty} \pi (g_n(a_2)) = \lim_{n \to \infty} \gamma_n (\pi(a_2)) = \lim_{n \to \infty} \gamma_n (x) = x^+_\infty,\]
which is a contradiction to (\ref{eq1}).\qed\\[\baselineskip]
\noindent{\bf Remark:} %Theorem \ref{notconical} seems independently interesting and apparently is previously unnoticed.
It is well-known \cite{BM}, that a finitely
generated Kleinian group without parabolic elements
is geometrically finite if and only if all limit points are conical
limit points. It is interesting that properties of the Cannon-Thurston
map apply to produce a natural class of
``explicit'' non-conical limit points for doubly degenerate surface
groups.\\[\baselineskip]
The key claims needed for the proof of Theorem \ref{surfacegroup} are now contained in the following proposition.

\begin{proposition}
\label{discrete}
\begin{enumerate}
\item ${\cal X}(\Gamma)$ is a locally finite collection of geodesics in ${\bf H}^3$.
\item If $\Gamma_1$ is a subgroup of finite index in
$\Gamma$ then ${\cal X}(\Gamma_1) = {\cal X}(\Gamma)$.
\item If $g\in \PSL(2,{\bf C})$, then  ${\cal X}(g\Gamma g^{-1}) = g{\cal X}(\Gamma)$.
\end{enumerate}
\end{proposition}
\noindent
{\bf Proof.} To prove (1), we suppose to the contrary that there exists a compact set $K \subset {\bf H}^3$ that intersects infinitely many geodesics in ${\cal X}(\Gamma)$ and arrive at a contradiction.
As ${\cal X}(\Gamma)$ is the orbit of finitely many geodesics $\{\xi^{(1)},\ldots,\xi^{(m)}\}$, there exists $1 \leq j \leq m$ and a sequence of distinct elements $\{\gamma_n\}_{n=1}^\infty \subset \Gamma$ so that $K \cap \gamma_n(\xi^{(j)}) \neq \emptyset$ for all $n$.  Equivalently, $\gamma_n^{-1}(K) \cap \xi^{(j)} \neq \emptyset$ for all $n$.  However, $\xi^{(j)}$ joins $x_+^{(j)}$ and $x_-^{(j)}$ and so by definition at least one of these points is a conical limit point.  This contradicts Lemma \ref{notconical}, and hence part (1) holds.\\[\baselineskip]
For part (2), we let $\Gamma_1 < \Gamma$ be a finite index subgroup and $p:S_1 \rightarrow S$ the associated finite sheeted cover.
The ending laminations for $\Gamma_1$ are the preimage $p^{-1}(\nu_\pm)$ of the ending laminations for $\Gamma$.
Since the universal coverings $q:\widetilde S \to S$ and $q_1:\widetilde S \to S_1$ satisfy $q = p \circ q_1$, it follows that $\nu_\pm$ and $p^{-1}(\nu_\pm)$ define the same set of principal regions in $\widetilde S$.  Since ${\cal X}(\Gamma_1)$ and ${\cal X}(\Gamma)$ are defined in terms of these regions, properties of the Cannon--Thurston map described in Theorem \ref{CTMapsforsurfaces} imply
${\cal X}(\Gamma_1) = {\cal X}(\Gamma)$. \\[\baselineskip]
To prove (3), we note that for a fixed $g \in \PSL(2,{\bf C})$, the equivariant map
associated to the conjugated  degenerate group $g \, \rho(\pi_1(S)) \, g^{-1}$
is obtained as a composition
$g \circ \iota : \widetilde S \longrightarrow {\bf H}^3$.
Indeed, for all $\gamma \in \pi_1(S)$, the action of $\gamma$ on $\widetilde S$ and on ${\bf H}^3$ via $\rho$ and its conjugate fit into the following commutative diagram
\[
\xymatrix@=1.5cm{ \, \widetilde S \, \ar[r]^{\large \iota} \ar[d]^\gamma & \, {\bf H}^3 \,  \ar[r]^g \ar[d]^{\rho(\gamma)} & \, {\bf H}^3 \, \ar[d]^{g.\rho(\gamma).g^{-1}}\\
\, \widetilde S \, \ar[r]^\iota & \, {\bf H}^3 \, \ar[r]^g & \, {\bf H}^3 \, \\}
\]

Now $g$ is a hyperbolic isometry and therefore admits an extension to ${\bf H}^3 \cup S^2_\infty$
and one sees that the two maps
$$ \widehat{g \circ \iota}, g \circ \hat{\iota} : \widetilde S \cup S^1_\infty \longrightarrow {\bf H}^3 \cup S^2_\infty$$
agree on the dense set $\widetilde S$, hence everywhere by continuity.  Since the geodesics in ${\cal X}$ are defined in terms of the boundary values of the Cannon--Thurston map, the result follows. \qed\\[\baselineskip]
The proof of the first part of Theorem \ref{surfacegroup} will follow from the
next claim. For if $C(\Gamma)$ is not discrete we can find a collection
of geodesics in ${\cal X}(\Gamma)$ violating (1) of  Proposition \ref{discrete}.\\[\baselineskip]
\noindent{\bf Claim 1:}~{\em Suppose that $g \in C(\Gamma)$, then $g$
leaves ${\cal X}(\Gamma)$ setwise invariant.}\\[\baselineskip]
\noindent{\bf Proof of Claim 1:}~Let $g\in C(\Gamma)$.
Then $\Gamma \cap g\Gamma g^{-1}$
has finite index in $\Gamma$ and $g\Gamma g^{-1}$. Proposition
\ref{discrete}(2) shows that the collection of special geodesics
associated to $\Gamma \cap g\Gamma g^{-1}$ is ${\cal X}(\Gamma)$, and
Proposition \ref{discrete} (3) shows that
$g({\cal X}(\Gamma))$ is the set of special geodesics associated
to $g\Gamma g^{-1}$. Proposition \ref{discrete}(2) applied to
$\Gamma \cap g\Gamma g^{-1}$ as a subgroup of $g\Gamma g^{-1}$ now completes
the proof of the claim.\qed\\

We now prove
\begin{theorem}
\label{fingen}
$C(\Gamma)$ is finitely generated.
\end{theorem}
\noindent
{\bf Proof:}  For $\xi \in {\cal X}(\Gamma)$, set $\Stab(\xi) = \{ g \in C(\Gamma) : g \xi = \xi \}.$
Since $C(\Gamma)$ is discrete, it follows that each $\Stab(\xi)$ is either a finite cyclic or dihedral group, or an extension of
an infinite cyclic or infinite dihedral group by a finite cyclic group.  In any case, $\Stab(\xi)$ is finitely generated.

Now recall that ${\cal X}(\Gamma)$ is the $\Gamma$ orbit of the finite set of geodesics $\{ \xi^{(1)},\ldots,\xi^{(m)}\}$.
Let ${\cal D}_r = \Stab(\xi^{(r)})$ for each $r = 1,\ldots,m$.
Then for any $\xi \in {\cal X}(\Gamma)$, there exists $\gamma \in \Gamma$ and $1 \leq r \leq m$ such that $\xi = \gamma(\xi^{(r)})$ and so
\[ \Stab(\xi) = \gamma  \Stab(\xi^{(r)}) \gamma^{-1} = \gamma {\cal D}_r \gamma^{-1}.\]

Let $K_1 = \langle {\cal D}_1, \ldots, {\cal D}_m, \Gamma \rangle $, and consider the action of $K_1$ on ${\cal X}(\Gamma)$.
Since $\Gamma < K_1$, this has at most $m$ orbits.
Suppose that there is an element $\alpha_1 \in C(\Gamma)$ which identifies
some of the orbits not identified by $K_1$. Set $K_2 = \langle \alpha_1 , K_1 \rangle$. We repeat this process: since
the number of orbits can only go down, we may continue in this way until reaching $K_j$, a subgroup
for which there are no extra orbit identifications (beyond those arising from $K_j$) possible for
{\em any} choice of an element of $C(\Gamma)$. The proof of the theorem is completed by the next Claim.\\[\baselineskip]
 \noindent{\bf Claim 2:}~ In this situation, $K_j = C(\Gamma)$ \\[\baselineskip]
 \noindent{\bf Proof of Claim 2:}~Pick $\beta \in C(\Gamma)$ and
choose any geodesic $\xi \in {\cal X}(\Gamma)$.
Then $\beta(\xi) = k(\xi)$ for
 some $k \in K_j$ so that $k^{-1}.\beta$ stabilises $\xi$, and
 therefore lies in $\gamma.{\cal D}_r.\gamma^{-1}$ for some $\gamma \in
 \Gamma$ and one of the stabiliser groups ${\cal D}_r$. Since
$\gamma.{\cal D}_r.\gamma^{-1} \leq K_j$, it follows that $\beta\in K_j$
so  proving the claim and
completing the proof of the theorem. \qed\\

\noindent {\bf Remark:}~ Theorem \ref{fingen} follows from the next theorem,
but it seemed interesting to give a direct proof using the structure of the special geodesics.\\
% \marginpar{\tiny The version of the covering theorem we're using
%   actually incorporates Agol's version.  This is the reason for the
%   rephrasing of this remark.  Also, local finiteness of ${\cal
%     X}(\Gamma)$ is used to deduce the structure of ${\cal D}_k$, in
%   particular its finite generation.}
%%%%%%%%%%%%%%%%%%%%%%%%%%%%%%%%%%%%%%%%%%%%%%
%
%
% Analysis of the case [C(\Gamma):\Gamma] = \infty
%
%
%%%%%%%%%%%%%%%%%%%%%%%%%%%%%%%%%%%%%%%%%%%%%%
We now determine when $[C(\Gamma):\Gamma]$ is finite, which will complete the proof of Theorem \ref{surfacegroup}.

\begin{theorem}
\label{infindexcomm}
The group $\Gamma$ is a fibre group if and only if $[C(\Gamma):\Gamma] = \infty$.
\end{theorem}
\noindent
{\bf Proof:}~Let $\Gamma_0 < \Gamma$ be the largest surface subgroup (which has index at most 2 inside $\Gamma$).
If $\Gamma$ is a fibre group then
$[C(\Gamma):\Gamma]$ is infinite. Indeed, $C(\Gamma_0) = C(\Gamma)$ contains
$N_{\PSL(2,{\bf C})}(\Gamma_0)$ (the normalizer of $\Gamma_0$ in $\PSL(2,{\bf C})$), which is a lattice, and we deduce that $C(\Gamma)$ is a
lattice.\\[\baselineskip]
Now suppose $[C(\Gamma):\Gamma] = \infty$, and consider the infinite sheeted covering
$$p: {\bf H}^3/\Gamma \longrightarrow    {\bf H}^3/C(\Gamma) $$
The manifold $  {\bf H}^3/\Gamma$ has at most two ends which are (both) geometrically infinite by assumption, and so $p$ is infinite-to-one on a neighborhood of at least one of these ends.
It follows from Theorem \ref{canary}, that $\Gamma$ is a fibre group.\qed

%%%%%%%%%%%%%%%%%%%%%%%%%%%%%%%%%%%%%%%%%%%
%%%%%%%%%%%%%%%%%%%%%%%%%%%%%%%%%%%%%%%%%%%
\section{Parabolics} \label{S:parabolics}
%%%%%%%%%%%%%%%%%%%%%%%%%%%%%%%%%%%%%%%%%%%
%%%%%%%%%%%%%%%%%%%%%%%%%%%%%%%%%%%%%%%%%%%

Here we explain the mild generalization of Theorem \ref{commensurator} where we allow our groups to contain parabolics.

\begin{theorem}
Suppose $\Gamma < \PSL(2,{\bf C})$ is a finitely generated torsion free Kleinian group of the first kind which is not a lattice.  If $\Gamma$ is not free and is not the fundamental group of a compression body with toroidal lower boundary, then $C(\Gamma)$ is discrete.  Moreover, $[C(\Gamma):\Gamma] = \infty$ if and only if $\Gamma$ is a fiber group.
\end{theorem}
\noindent
{\bf Proof:}~(Sketch)  If $\Gamma$ has no parabolics, then this reduces to Theorem \ref{commensurator}, so we assume that $\Gamma$ contains parabolic elements.

We refer the reader to \cite{Ca0} and \cite{Ca1} for terminology and a detailed discussion of the notation used here.  Let $q:{\bf H}^3 \to N = {\bf H}^3/\Gamma$ be the universal covering, $\epsilon > 0$ some number less than the 3--dimensional Margulis constant, and let $N^0_\epsilon$ denote the result of removing the cuspidal $\epsilon$--thin part of $N$; the Tameness Theorem implies $N^0_\epsilon$ is tame (\cite{A},\cite{CG}).  Let $M \subset N^0_\epsilon$ denote a relative compact core for $N^0_\epsilon$: a compact core for which $P = M \cap \partial N^0_\epsilon$ is a disjoint union of annuli and tori in $\partial M$ and so that the ends of $N^0_\epsilon$ are in a one-to-one correspondence with the components of $\partial M - P$.

For any group $\Gamma$ as in the statement of the theorem, we need to find a finite set ${\cal P}_N$ of finite area pleated surfaces ${\cal P}_N$ in $N$, such that the locally finite collection of pleated surfaces
\[ {\cal P}_\Gamma = q^{-1}({\cal  P}_N) \]
satisfies the following two properties (c.f.~(\ref{conj eqn1}) and (\ref{fi eqn1}) from \S \ref{S:not surface}).
\begin{enumerate}
\item For any finite index subgroup $\hat \Gamma < \Gamma$ we have ${\cal P}_{\hat \Gamma} = {\cal P}_\Gamma$,
\item For all $g \in \PSL(2,{\bf C})$, we have ${\cal P}_{g \Gamma g^{-1}} = g ({\cal P}_\Gamma)$, and
\end{enumerate}
Given these two properties, the proof of the Theorem follows as in the proof of Theorem \ref{notsurface}.

First, suppose that $\Gamma$ is not virtually a surface group and fix a nontoroidal component $S \subset \partial M$, which exists since $\Gamma$ is not a lattice.
Let $B_S$ be the compression body associated to $S$ from \cite{Bon0} as in the proof of Theorem \ref{notsurface}.
By hypothesis, $B_S$ is not a handlebody and so has some nontrivial lower boundary $\partial_-B_S \neq \emptyset$.
Furthermore, we claim that $\partial_- B_S$ cannot be a union of tori.
To see this, note that this would imply that either we could take $B_S = M$, or else some component $T \subset \partial_- B_S$ is not peripheral in $M$.
The former case is ruled out by the hypothesis that $\Gamma = \pi_1(M)$ is not the fundamental group of a compression body with toroidal lower boundary, while the latter case is impossible since $M$ is the core of a hyperbolic manifold, so has no nonperipheral incompressible tori.

Let $\Sigma_1,\ldots,\Sigma_k \subset \partial_- B_S$ denote the set of non-toroidal components (which is nonempty) and let $\Gamma(1,S),\ldots,\Gamma(k,S) < \Gamma$ be the associated nonempty, finite collection of surface subgroups.
Since we have assumed that $\Gamma$ is not virtually a surface group, we have $[\Gamma:\Gamma(j,S)] = \infty$ for all $1 \leq j \leq k$.

For each $1 \leq j \leq k$, let $N(j,S) = {\bf H}^3/\Gamma(j,S)$ and $N_\epsilon^0(j,S)$ the complement of the cuspidal $\epsilon$--thin part in $N(j,S)$.
Appealing to the covering theorem (for $3$--manifolds with cusps; see \cite{Ca1,CL}), we again deduce that at least one end of $N_\epsilon^0(j,S)$ is geometrically finite.
The boundary of the convex core $C_{\Gamma(j,S)}$ maps into $N$, and the set of all such surfaces is a finite collection of finite area hyperbolic surfaces canonically associated to $S$.
Let ${\cal P}_N$ be the set of these pleated surfaces over all components $S \subset \partial M$.  Properties 1 and 2 are easily verified for this family, and discreteness of $C(\Gamma)$ follows.

If $\Gamma$ is virtually a surface group, then since we have assumed $\Gamma$ has parabolics, $P \subset \partial M$ is a nonempty collection of annuli.  For simplicity, assume $M = S \times [0,1]$ (if $M$ is a twisted $I$--bundle, take the two-fold cover).
Let $\Sigma \subset \partial M - P$ be a component whose closure nontrivially meets $P$ and let $\Gamma_\Sigma \subset \Gamma$ denote the (injective) image of $\pi_1(\Sigma)$ in $\pi_1(S) = \Gamma$.
There is at least one such $\Sigma$ since $P \neq \emptyset$.
Since $\Sigma$ corresponds to a proper subsurface of $S$, $\Gamma_\Sigma$ has infinite index in  $\Gamma$.
Letting $N(\Sigma) = {\bf H}^3/\Gamma_\Sigma$ and $N^0_\epsilon(\Sigma)$ be the complement of the $\epsilon$--thin cuspidal part, arguing as above one sees that $N^0_\epsilon(\Sigma)$ has at least one geometrically finite end.
Therefore, we can associate to $\Sigma$ the components of $C_{\Gamma_\Sigma}$ mapped into $N$, which is a finite set of finite area pleated surfaces.
The set of all of these taken over all components $\Sigma \subset \partial M - P$ as above gives the required set of pleated surfaces.\\

Setting ${\cal P}$ to be the union of all the surfaces in ${\cal P}_\Gamma$ we have
\[[C(\Gamma):\Gamma] = \frac{\Area( {\cal P}/\Gamma)}{\Area({\cal P}/C(\Gamma))}.\]
It follows that under the hypotheses of the theorem, if $\Gamma$ has parabolics, $[C(\Gamma):\Gamma]< \infty$.
If $\Gamma$ has no parabolics, then applying Theorem \ref{commensurator} completes the proof. \qed

%%%%%%%%%%%%%%%%%%%%%%%%%%%%%%%%%%%%%%%%%%%%%%%%%%%%%
%
%
% Bibliography
%
%
%%%%%%%%%%%%%%%%%%%%%%%%%%%%%%%%%%%%%%%%%%%%%%%%%%%%%

\bigskip

 \noindent
Department of Mathematics,\\ University of Illinois at Urbana-Champaign,\\
Urbana IL 61801

\noindent Email:~clein@math.uiuc.edu\\[\baselineskip]
 \noindent
 Department of Mathematics,\\ University of California,\\ Santa Barbara, CA
93106.

\noindent Email:~long@math.ucsb.edu\\[\baselineskip]
 Department of Mathematics,\\
 University of Texas,\\
 Austin, TX 78712.

\noindent Email:~areid@math.utexas.edu\\

\end{document}